\documentstyle{amsppt}
\input amstex
\magnification = 1000
\baselineskip = 16pt
\pagewidth{6.5  true in}
\pageheight{8.5 true in}
\overfullrule=0pt
\parskip = 3pt
\addto \tenpoint {\normalbaselineskip18pt\normalbaselines}
\pageno = 1
\TagsOnRight
\NoBlackBoxes

\document
\topmatter
\title  ON THE  FIBER   CONES  OF  GOOD  FILTRATIONS  \endtitle                        
    
\author  Duong Quoc Viet  \endauthor
\affil{ Department of Mathematics, Hanoi University of Education}\\ 
{136-Xuan Thuy Street, Hanoi, Vietnam}\\ 
{\sevenrm E-mail: duongquocviet\@\ fmail.vnn.vn }\\
\endaffil
\endtopmatter
{\leftskip = 1.5cm \rightskip = 2cm
{{\eightrm 
\noindent
{\bf\sevenrm ABSTRACT }:
Let $A$ be a Noetherian local ring with the maximal ideal ${\frak m}$ 
and  an ${\frak m}$-primary ideal $J$.  Let $\Cal  F = \{I_n\}_{n \ge 0}$  be a good  filtration  of ideals in $A.$ 
Denote by 
$F_J(\Cal F) = \underset{n\ge 0} \to \bigoplus (I_n/JI_n)t^n$ 
the  fiber 
cone of $\Cal  F$  with respect to
$J.$  
The paper characterizes the multiplicity   and 
the Cohen-Macaulayness of $F_J(\Cal  F)$ in terms of minimal reductions of  $\Cal  F.$ 
} \par }}

\document
\vskip 0.3cm
 
{\bf \heading 1. Introduction \endheading \smallskip}
\vskip 0.2cm
\pagewidth{6.5 true in}
\noindent
Let $(A, \frak m)$ be a Noetherian local ring 
with  maximal ideal $\frak m,$  infinite residue field  $k = A/\frak m,$  $\Cal  F = \{I_n\}_{n \ge 0}$  a good   filtration of ideals in $A.$  
  Let $J$ be an  arbitrary ${\frak m}$-primary ideal.  
 Define 
$F_J(\Cal F) = \underset{n\ge 0} \to \bigoplus (I_n/JI_n)t^n$ to be the {\it fiber 
cone } of $\Cal F$  with respect to
$J.$ In the case where  $\Cal F = \{I^n\}_{n \ge 0}$  is an $I$-adic filtraion,  $F_J(\Cal F)$ is  called the {\it fiber 
cone } of $I$  with respect to
$J$  and it  is  denoted by $F_J(I).$      
The notations $\ell(I) = \dim F_{\frak m}(I)$ and $\ell(\Cal F) = \dim F_{\frak m}(\Cal F)$ will  mean  the 
{\it analytic spread } of $I$ and of  $\Cal F,$ respectively. 

\footnotetext[] 
{\leftskip = 0.5 cm \rightskip = 1cm
\item{ }{\it Mathematics  Subject  Classification} (2000): Primary 13H10. Secondary 13A15, 13A30, 13C14, 13H15. 
\item{ }$ Key\; words \;  and \; phrases:$ Cohen-Macaulay ring,  Multiplicity, Fiber  cone,  Filtration, weak-(FC)-sequence.}
The multiplicity and  the Cohen-Macaulayness  of  fiber cones   
are  usually  interesting problems.   
These  problems  are  concerned by many  authors 
in the past years.  Using different approaches, the authors  investigated the  Cohen-Macaulayness and other properties of  fiber cones $F_{\frak m}(I)$ and $F_J(\Cal F)$  (see for instance [C-G-P-U], [C-P-V], [C-Z],  [Co-Z],  [D-R-V],  [H-H],  [J-V1],  [J-P-V], [C],  [Vi2]).      
Using weak-(FC)-sequences of  ideals in  local rings, the author of  [Vi2] characterized the  multiplicity and  the Cohen-Macaulayness of $F_{\frak m}(I)$ in terms of  minimal reductions of ideals.
 The results of [Vi2] recovered  some  earlier  results  of  Huneke-Sally [H-Sa],  Shah [Sh1, Sh2] and Cortadellas-Zarzuela [C-Z].  

In terms of  minimal  reductions of  filtrations, the aim of this  paper  is to  give characterizations  of the multiplicity (Theorem  3.3, Section 3)  and  the  Cohen-Macaulayness (Theorem  4.2, Section 4) of  the fiber cone $F_J(\Cal F)$  of a  good  filtration $\Cal F.$  
  A crucial role in this paper is played  by the use of  weak-(FC)-sequences  of  good  filtrations (see Section 2). The results of this paper prove  that the main results of  [Vi2] are  still  true for  fiber cones of  good  filtrations.  Moreover, from  the  main  result  we show that: For  any  good  filtration $\Cal  F = \{I_n\}_{n \ge 0}$ with $\ell(\Cal F) =1,$ $\underset{n\ge 0} \to \bigoplus (I_{Ln}/\frak mI_{Ln})t^n$            
 is Cohen-Macaulay for  all large  $L$ (Corollary 4.3,  Section 4); and we obtain  more  favorite
results  than the results in  [Vi2] (Remark  4.4).

This paper is divided into  four  sections. Sect.2 deals  with   weak-(FC)-sequences of good filtrations.
Sect.3 investigates  the multiplicity of fiber cones of  good filtrations.  Sect.4  is devoted to the discussion of
the Cohen-Macaulayness of fiber cones of  good  filtrations.

\vskip 0.3cm 

{\bf \heading{2. The weak-(FC)-sequences of good  filtrations }
\endheading \smallskip}

\vskip 0.2cm  
\noindent
 The author in [Vi1] built (FC)-sequences of ideals in local rings for calculating mixed multiplicities. In
order to study the multiplicity and the Cohen-Macaulayness for   fiber  cones of good   filtrations, this section
introduces  weak-(FC)-sequences of good   filtrations and some important properties of these sequences.

A  filtration $\Cal  F = \{I_n\}_{n \ge 0}$  of ideals in $A$ is   a   chain  of  ideals  $I_n$ such that 
$A = I_0,$ $I_1$ is   a  proper  ideal  of  $A,$  $I_{n+1} \subseteq I_n $ and $I_nI_m \subseteq I_{n+m}$ for  all  $n, m \ge 0.$ Let $I$ be 
an  ideal of $A.$  $\Cal F$ is called
 an  {\it $I$-good  filtration}  if  $II_n \subseteq I_{n+1}$  for  all  $n \ge 0$ and  $I_{n+1} = II_n $  for   all  large $n.$ In  this  case, $I \subseteq I_1.$   $\Cal F$ is called  a {good  filtration} if it  is an $I$-good  filtration   for  some  ideal  $I$  of $A.$ It  is  easily seen that  $\Cal F$ is   a good  filtration if and only if $\Cal F$ is  an  $I_1$-good  filtration.  A  good  filtration  $\Cal  F$ is called  a  nilpotent  filtration  if  $I_n = 0$ for  all  large $n.$ This is  equivalent  to  $I_1$ nilpotent.  Set $\Cal F/I = \{I_n(A/I)\}_{n \ge 0}$ and  $ 
F_J(\Cal F/I ) = \underset{n\ge 0} \to \bigoplus \big[(I_n+I)/(JI_n+I)\big]t^n$  for  any ideal  $I$ of $A.$  
 
{\proclaim{{\bf Definition }} 
Let $I$ be  an ideal of $A$  and $\Cal  F = \{I_n\}_{n \ge 0}$  a non-nilpotent good  filtration  of ideals of $A.$   
An element $x \in  I_1$  is called 
a weak-(FC)-element  with respect to
$( I, \Cal F) $ 
if the following conditions are satisfied: 

\noindent
\item{ }{\rm(FC$_1$):}   
$I^m I_n \bigcap (x) = 
I^m x I_{n-1}$ for  all   large $n$ and for all non-negative integers $m.$

\noindent
\item{ }{\rm(FC$_2$):}  $x$ is  a filter-regular  element with respect to $I_1,$ 
 i.e., 
$ 0: x \subseteq 0:I_1^{\infty}.$ 

\noindent
Let  $x_1, x_2, \ldots,x_s$ be a sequence in  
  $I_1.$ 
For each $i = 0,1, \ldots,s,$ 
set
\noindent 
$A_i = A/(x_1,x_2, \ldots,x_i);$ 
${\bar x}_{i+1}$ the image  of $x_{i+1}$ in  ${A_i};$ 
 $ \Cal F_i   = \Cal F/(x_1,x_2, \ldots,x_i).$   Then

\item{{\rm(i)}} The  sequence $x_1,\ldots,x_s$ is called  a 
weak-(FC)-sequence 
in $I_1$ with respect to 
$(I, \Cal F) $ 
 if 
${\bar x}_{i+1}$ is a 
weak-(FC)-element  with respect to
 $(IA_i, \Cal F_i)$                             
for each $i = 0,1, \ldots,s-1.$

\item{{\rm(ii)}} A  weak-(FC)-sequence $x_1,\ldots,x_s$ in $I_1$ 
with respect to $(I, \Cal F) $ 
is called a  maximal weak-(FC)-sequence  if 
$\Cal F_{s-1}$ is  a non-nilpotent   filtration   of $A_{s-1},$ but 
$\Cal F_s$ is  a  nilpotent   filtration  of $A_s.$  
\endproclaim

An  ideal $\frak I$ of $A$ is  called  a {\it reduction} of   a  good  filtration $\Cal F = \{I_n\}_{n \ge 0}$  if  $\Cal F$ 
is  an  $\frak I$-good  filtration.    
The least integer $n$ such that ${\frak I} I_n =  I_{n+1}$ is called  the {\it  reduction number  } of $\Cal F$ with respect to $\frak I$ and we denote this integer  by $r_{\frak I} (\Cal F).$ A  reduction $\frak  I$ of $\Cal F$ is called  a {\it  minimal reduction } if it does not properly contain  any other reduction of $\Cal F.$ The {\it reduction number } of $\Cal F$ is defined by 
$$r(\Cal F) = \min\{ r_{\frak I }(\Cal F) \;  \mid \; {\frak I}   \text { is a minimal reduction of } {\Cal F}  \}.$$
 Note that in the case of $\Cal F$ an  $I$-adic filtration, $\frak I$  is  called  a {\it reduction} of $I$ and the notations  $r_{\frak I} (I)$ and $r(I)$  will  mean  the   reduction number  of $I$ with respect to $\frak I$  and  the   reduction number  of $I$, respectively.
Northcott  and  Rees in [N-R]  showed  that  a reduction  $\frak  I$  of  $I$ is  a  minimal reduction if  and  only if 
 the minimal number  of  generators  of  $\frak I$ 
is  equal to  the  analytic spread  $\ell(I) = \dim  F_{\frak m}(I)$  of  $I.$  
If ht $I$ is the height of $I$ then 
$\text{\rm ht }(I) \leqslant \ell(I) .$ In the case of 
$\text{\rm ht }(I) = \ell(I),$ $I$ is called {\it equimultiple}.  A  good  filtration $\Cal F$  is called  an {\it equimultiple  filtration} if $\text{\rm ht }(I_1) = \ell(\Cal F) .$ 

Define $R(\Cal F) = \bigoplus_{n\ge 0}I_nt^n$ and $G(\Cal F) =  \underset{n\ge0}\to\bigoplus (I_n/I_{n+1})t^n$   to  be the  Rees  algebra and  
the associated graded ring of $\Cal F,$ respectively. Denote by $\frak M$ the maximal homogeneous ideal  of $R(\Cal F).$

Now, we  briefly  give some comments  on  weak-(FC)-sequences  of  a good  filtration of  ideals in  $A$ and  the fiber cone of  good filtrations  by   the  following  remark.

\noindent
{\bf Remark 2.1.} 
\noindent
\item{ }{\rm (i)} 
If  $\Cal F$ is a nilpotent good  filtration  of ideals of $A$  then  $I_1^n \subseteq I_n = 0$ for all large $n$.   Consequently, 
  for  any  element $x\in I_1$ and  for  any ideal $I$ of $A,$  we always have   
$I^m I_n \bigcap (x) = 0 = 
{I}^m x I_{n-1}$ for all large $n$  and  $0: x  \subset  A   =  0: {I_1}^{\infty}.$  Hence  the conditions \rm(FC$_1$) and \rm(FC$_2$) always  are satisfied  for  all $x\in I_1.$    
    This only obstructs and does not carry useful. That  is  why  in  definition  of weak-(FC)-elements,  we have  to  exclude  the case where  
$\Cal F$ is a nilpotent good  filtration  of ideals of $A.$
\noindent
\item{ }{\rm (ii)} Since  $\Cal F$ is a  good  filtration  of ideals of $A,$  there  exists  $u$ such that $I_{n} = I_1^{n-u}I_{u}$  for  all  large $n.$    
 By Artin-Rees lemma, there exists integer $v$ such that
$$(0: {I_1}^{\infty}) \bigcap I_n  \subseteq  (0: {I_1}^{\infty}) \bigcap  I_1^{n-u} = ((0: {I_1}^{\infty}) \bigcap I_1^v) I_1^{n-u-v} \subseteq (0: {I_1}^{\infty})I_1^{n-u-v} $$ for all $n-u  \ge v.$ 
Since $I_1^{n-u-v}(0: {I_1}^{\infty}) = 0$ for all large 
enough 
$n,$
$(0: {I_1}^{\infty}) \bigcap I_n = 0$ for all large  $n.$  
\noindent
\item{ }{\rm (iii)}  Suppose that $x\in I_1$ is  a filter-regular element with respect to $I_1.$   Consider $$\lambda_x: I_n \longrightarrow xI_n, \;  y\mapsto xy.$$ It is clear that $\lambda_x$ is  surjective and $\text {Ker } \lambda_x =(0 : x) \bigcap I_n$.  Since $x$ is a filter-regular element,  $$\text {Ker } \lambda_x = (0 : x) \bigcap I_n  \subseteq (0: {I_1}^{\infty}) \bigcap I_n =  0 $$ for all large $n$. Therefore, $xI_n \simeq I_n.$ This  follows that  $xII_n \simeq II_n$ for  any ideal $I$ of $A$  and  for all large $n$.  Hence   for  any  ideal $I$ of $A,$ we have an isomorphism  of  $A$-modules 
$I_n/II_n \simeq xI_n/xII_n\;\;\text{ for  all large }n. $

\noindent
\item{ }{\rm (iv)}  Set
$ A^* = A/ 0:I_1^{\infty};$ $ I^* =  IA^*;$   
$I_n^* = I_nA^*$  for  all  $n \ge 0;$  
 $a^*$ the image 
of  $a \in A$ in  $A^*.$  Suppose that $x\in I_1$ is  a filter-regular element with respect to $I_1.$ Since $0 : x \subseteq 0:I_1^{\infty}$, $x^*$ is  a non-zero-divisor in $A^*.$ Hence   
$I_n^*/I^*I_n^* \simeq x^*I_n^*/x^*I^*I_n^*\;\;\text{ for  all }n. $

\noindent
\item{ }{\rm (v)} If  $\ell(\Cal F) = 1$ and  an  element $x$ such that $(x)$  is a reduction of $\Cal F,$ then for  any  ideal $I$ of $A$  and  for  all   large $n,$  we have  $ I^m I_n \bigcap (x)  = I^m xI_{n-1}\bigcap (x) = I^m xI_{n-1}$ for all non-negative integers $m.$ 
On  the other  hand  $0 : x   \subseteq 0: {I_1}^{\infty}.$ Hence  $x$ is  
a weak-(FC)-element  with respect to
$( I, \Cal F). $ 

\noindent
\item{ }{\rm (vi)} It is easily seen that if 
$J$ is  an ${\frak m}$-primary ideal  of  $A,$ then in $R(\Cal F)$ we  have  $\sqrt{JR(\Cal F)}  = {\frak m }R(\Cal F).$
Hence $\dim F_J(\Cal F) = \dim F_{\frak m}(\Cal F) = \ell(I_1) = \dim  F_J(I_1),$ and 
  if $\frak I$ is a minimal reduction of $\Cal F,$ then $\frak It F_J(\Cal F)_{\frak M}$ is  an  ideal  of  parameter  for  $F_J(\Cal F)_{\frak M}.$  

\noindent
\item{ }{\rm (vii)} Let  $J$ be  an ${\frak m}$-primary ideal  of  $A.$  
Set  $S_j = \bigoplus_{n\ge j}(I_n/JI_n)t^n.$ Then $S_j$ has  a  natural $F_J(I_1)$-module structure  given  by     
$(a + JI_1^m)(x +JI_n) = (ax + JI_{m+n})$ for $a \in I_1^m , \;  x \in I_n.$   
Since $\Cal F$ is a  good  filtration  of ideals of $A,$  it is easily seen that  
$S_j = \bigoplus_{n\ge j}(I_n/JI_n)t^n  =  F_J(I_1)(I_j/JI_j)t^j$ for   all  large $j.$   Hence $l_A(I_n/JI_n)$    
is a polynomial $Q(n)$  for   all  large $n,$  and $\deg Q(n) = \dim F_J(\Cal F) -1 = \ell (I_1) -1.$
 
\noindent
\item{ }{\rm (viii)}
Let $I$ be  an ideal of $A$  and $\Cal  F = \{I_n\}_{n \ge 0}$  a  good  filtration  of ideals of $A.$  Recall that an element $a \in  I_1$  is called 
a {\it superficial element } with respect to
$( I, \Cal F) $ if there exists a
non-negative integer $c$  such that
$(I^mI_{n+1}:  a) \bigcap I^mI_c  = I^mI_n$ 
for  all  $n \ge c$ and  for all non-negative integers $m.$ 
The notion  of the superficial  elements goes  back  to  P. Samuel  [Z-S].  The classical theory of the superficial elements  becomes  an  important tool in  local  algebra and has been continually developed (see for  instance  [H-S], [K], [R-S]).  
We now show that if an element $x\in I_1$ is a weak-(FC)-element with respect to $(I,\Cal F)$ then $x$ is also a superficial element with respect to $(I, \Cal F)$. 
Indeed, if $x\in I_1$ is a weak-(FC)-element with respect to $(I, \Cal F)$, then for all large $n$ and all $m\geqslant 0$, we have
$$I^mI_n:x =\big(I^mI_n\cap (x)\big):x =xI^mI_{n-1}:x
=I^mI_{n-1} +(0:x)
\subseteq I^mI_{n-1} + (0:I_1^\infty).
$$
By (ii) , there exists a positive integer $c$ such that $(0:I_1^\infty)\cap I_c=0.$ Thus for all large $n\ (n>c)$ and all $m\geqslant 0$,
$$(I^mI_{n}:x)\cap I^mI_c \subseteq [I^mI_{n-1} + (0:I_1^\infty)]\cap I^mI_c
=I^mI_{n-1} + (0:I_1^\infty)\cap I^mI_c
=I^mI_{n-1}.
$$
The reverse inclusions are trivial. Hence $x$ is a superficial element with respect to $(I,\Cal{F})$.

A  minor variation in Rees's  argument [see the proof of Lemma 1.2, Re]  yields  the following lemma. 

\noindent
\proclaim{\bf Lemma 2.2(Generalized Rees's Lemma)} Let  $\Cal F = \{I_n\}_{n \ge 0}$  be a good  filtration of ideals of $A.$ 
Let $\frak I$ be  a reduction of  $\Cal F,$ $I$  an ideal of $A.$   
Let $\Sigma$ be a finite collection of prime ideals of $A$ not containing $I_1.$
Then  there exists an element $x \in {\frak I}\backslash \bigcup_{P \in \Sigma}P  $
 such that
$(x)\bigcap I^{m}I_n 
= I^{m}xI_{n-1}$ 
for  all   large  n  and   all  non-negative  integers  $ m. $
\endproclaim
\vskip0.2cm

Using Lemma 2.2,  we  will  show  that  the existence of weak-(FC)-sequences in  good  filtrations  is universal by  the following proposition.

\proclaim{\bf Proposition 2.3} Let $I$ be an ideal of $A,$ $\frak I$  a reduction of  $\Cal F = \{I_n\}_{n \ge 0}.$  Suppose that $\Cal F$  is non-nilpotent.
Then  there exists a weak-(FC)-element  in $ \frak I$ with respect to $(I, \Cal F).$
\endproclaim
 
\vskip0.2cm
\demo{Proof} Set  $\Sigma = \text{Ass}_A[A/ (0: I_1^{\infty})].$  It  is  easily seen that  
 $\Sigma =  \{P \in \text{Ass}A \mid P \nsupseteq  I_1\}$ and  
   $\Sigma$ is  finite. Since $I_1$ is non-nilpotent,  there exists 
$x \in \frak I $ such that $x \notin P$ for all $P \in \Sigma$ and  
$(x)\bigcap I^mI_n 
= I^mxI_{n-1}$  for all large $n$ and all non-negative integers $m$ by Lemma 2.2. 
Thus,   $x$ satisfies  condiction (FC$_1$).
Since $x \notin P$ for all $P \in \Sigma$, 
 $0 : x \subset 0 :{I_1}^\infty.$ Hence  $x$ satisfies  condition (FC$_2$). $\square$ 
\enddemo

Let $J$ be an ${\frak m}$-primary ideal of  $A.$
Denote by $\frak M$ the maximal homogeneous ideal  of  the Rees  algebra $R(\Cal F)$ of  a  filtration $\Cal F.$ 
The notation $e(F_J(\Cal F))$ will  mean  the  Hilbert-Samuel   multiplicity  of  local ring 
$F_J(\Cal F)_{\frak M},$ and it  is  called  the multiplicity of  $F_J(\Cal F).$   

\proclaim{\bf Proposition 2.4}  Let $J$ be an ${\frak m}$-primary ideal of  $A$. Let  $\Cal F = \{I_n\}_{n \ge 0}$  be a good  filtration of ideals in $A.$  Set $\ell = \ell(\Cal F)$  and  $F_J(\Cal F)^+ = \underset{n > 0} \to \bigoplus (I_n/JI_n)t^n.$  Then we have: 

\noindent
\item{ }{\rm(i)} If  $\ell  > 0$ then  $e(F_J(\Cal F)) = e(F_J(\Cal F)^+;\;  F_J(\Cal F)) = \lim_{n\to +\infty}\dfrac{(\ell -1)! \big[l_A(I_n/JI_n)\big]}{n^{\ell-1}}. $

\noindent
\item{ }{\rm(ii)} If $\ell  > 1$ and $x\in I_1$ is a weak-(FC)-element  with 
respect to
$( J, \Cal F)$  then $e(F_J(\Cal F/(x))) = e(F_J(\Cal F)).$ 

\noindent
\item{ }{\rm(iii)}  The length of maximal weak-(FC)-sequences in $I_1$  with respect to $(J, \Cal F)$  is \; $\ell.$   

\noindent
\item{ }{\rm(iv)}  If $x_1,x_2,\ldots, x_{\ell}$ is  a  
weak-(FC)-sequence in $I_1$ 
with respect to $( J,  \Cal F)$  then $(x_1,x_2,\ldots, x_{\ell})$  is  a  
\item{ }{\;\;\;\;\;}  minimal reduction of   $ \Cal F.$ 
\noindent
\item{ }{\rm(v)}  Any minimal reduction of $\Cal F $ is generated by a maximal  weak-(FC)-sequence in $I_1$  with respect to    
\item{ }{\;\;\;\;\;}$(J, \Cal F).$
\endproclaim
\vskip0.2cm
\demo{Proof} Since $J$ is  an ${\frak m}$-primary ideal of  $A,$ ${\frak m}/J$ is  a  nilpotent  ideal  in  $A/J$ , and  hence 
 $F_J(\Cal F)^+$ is a reduction  of $(\frak m/J)\bigoplus F_J(\Cal F)^+.$  This  gives  $e(F_J(\Cal F)) = e((\frak m/J)\bigoplus F_J(\Cal F)^+; F_J(\Cal F)) = e(F_J(\Cal F)^+; F_J(\Cal F)).$ By Remark 2.1 (vii),
there exists  a positive  integer  $j$ such that $S_j = \bigoplus_{n\ge j}(I_n/JI_n)t^n  =  F_J(I_1)(I_j/JI_j)t^j.$ 
Concider the exact sequence $$0\to S_j \to  F_J(\Cal F) \to \bigoplus_{0 \leqslant n \leqslant j-1}(I_n/JI_n)t^n \to 0.$$ Since 
$\dim \bigoplus_{0 \leqslant n \leqslant j-1}(I_n/JI_n)t^n = 0  < \ell =  \dim  F_J(\Cal F),$ it  follows  that $e(F_J(\Cal F)) = e(S_j).$ Direct computation shows that  $$ l_A\big[S_j/( F_J(\Cal F)^+)^{n+1}S_j \big] - l_A\big[S_j/( F_J(\Cal F)^+)^n S_j  \big] = l_A\big[(I_{n+j}/JI_{n+j})\big]$$ for  all   large enough $n.$ 
 Remember  that $l_A\big[(I_{n+j}/JI_{n+j})\big]$  is  a polynomial  of  degree $(\ell-1)$ for  all   large enough $n.$ Hence  
we get     
$e(F_J(\Cal F)) = e(S_j) = \lim_{n\to +\infty}\dfrac{(\ell -1)! \big[l_A(I_n/JI_n)\big]}{n^{\ell-1}}.  $ 
 This completes the proof of (i).  Let  $x\in I_1$ be a weak-(FC)-element  with 
respect to
$( J, \Cal F).$ 
Set 
$E_n = I_n(A/xA)$  for  all $n \ge 0.$ 
Then $\Cal F/(x) = \{E_n\}_{n\ge 0}$ and  for all large enough $n,$  
we have
$$\aligned l_A({{E_n}\over{JE_n}})  
 = l_A({{I_n}+(x)\over{ JI_n +( x)}})
&= l_{A}({{I_n}\over{JI_n +(x)\cap I_n}})\\ 
 &= l_A({{I_n}\over{JI_n}}) - l_A({{JI_n+(x)\cap I_n}\over{ JI_n}})\\
&= l_{A}({{I_n}\over{JI_n}})- 
  l_A({{(x)\cap I_n}\over{ JI_n \bigcap (x)}})\\
&= l_{A}({{I_n}\over{JI_n}})- 
  l_A({{xI_{n-1}}\over{x JI_{n-1}}})  \text{(because }\; x \; \text {satisfies the condition }  \text{(FC}_1\text{))}\\ 
&= l_{A}({{I_n}\over{JI_n}})- 
  l_A({{I_{n-1}}\over{ JI_{n-1}}}) \text{(by  Remark 2.1 (iii))}.\\ 
\endaligned $$
Consequently, it  holds  that   $$ l_A({{E_n}\over{JE_n}}) = l_{A}({{I_n}\over{JI_n}})- 
  l_A({{I_{n-1}}\over{ JI_{n-1}}}) \eqno (1) $$  for all 
large enough  $n.$ By (1) we get    
$e(F_J(\Cal F/(x))) = e(F_J(\Cal F)).$ Therefore,  we have(ii).  It  follows  readily  from (1)  that  $$\dim F_J(\Cal F/(x)) = \dim F_J(\Cal F)-1 \text{(**)}.$$  Next we prove (iii):  First note that if $\ell = 0$ then $I_1$     
is a nilpotent ideal of  $A.$  In this case,  the length of maximal weak-(FC)-sequences in $I_1$  with respect to $(J,\Cal F)$  is  $0 = \ell.$
If $\ell > 0$ and  $x_1,x_2,\ldots, x_s$ is a  
weak-(FC)-sequence in $I_1$ 
with respect to $( J,  \Cal F),$ then  using  (**),  we  easily  prove  by induction on $s$  that $$\dim F_J(\Cal F/(x_1,x_2,\ldots, x_s)) = \dim F_J(\Cal F)-s.$$  This  implies  that  the length of maximal weak-(FC)-sequences in $I_1$  with respect to $(J, \Cal F)$  is  $\ell.$ We get (iii). The proof of (iv): Let $x_1,x_2,\ldots, x_{\ell}$ be  a  
weak-(FC)-sequence in $I_1$ 
with respect to $( J,  \Cal F).$ The proof is by induction on $i \leqslant \ell $ that $(x_1,x_2,\ldots, x_{i}) \bigcap I_{n+1} = (x_1,x_2,\ldots, x_{i}) I_n$ for all large $n.$ The  case of $i=0$ is trivial.  Denote by $x_i'$ the image of $x_i$ in $A' =A/(x_1,x_2,\ldots, x_{i-1}).$ For suppose the result has been proved for $i-1 \ge 0.$ Set $N = (x_1,x_2,\ldots, x_{i-1}).$  Since $x_i'$  is  a weak-(FC)-element 
with respecto  $( JA',  \Cal F/N),$
$(x_i')\bigcap I_{n+1}A' =  x_i'I_{n}A'$ for  all  large $n.$  So 
$[I_{n+1}+N]\bigcap [(x_i) +N] = x_iI_n + N.$ Hence $$I_{n+1}\bigcap [I_{n+1}+N]\bigcap [(x_i) +N] = I_{n+1}\bigcap (x_iI_n + N).$$ This is equivalent to          
$I_{n+1}\bigcap [(x_i) +N] = x_iI_n + N \bigcap I_{n+1}.$ By inductive assumption, $N \bigcap I_{n+1}= NI_n$ for all large $n.$ 
Hence $(x_1,x_2,\ldots, x_{i})\bigcap I_{n+1} = (x_1,x_2,\ldots, x_{i})I_n$ for all large $n.$ The induction is complete. This gives   
$(x_1,x_2,\ldots, x_{\ell})\bigcap I_{n+1} = (x_1,x_2,\ldots, x_{\ell})I_n$ for all large $n.$
Since  $x_1,x_2,\ldots, x_{\ell}$  is  a  maximal weak-(FC)-sequence in $I_1$  with respect to $(J, \Cal F)$ by (iii),   
 $\Cal F/ (x_1,x_2,\ldots, x_{\ell})$ is a  nilpotent  filtration  of  $A/(x_1,x_2,\ldots, x_{\ell}).$ Therefore,   
$I_n \subseteq  (x_1,x_2,\ldots, x_{\ell})$  for  all 
large $n.$ Consequently, $$ I_{n+1} = (x_1,x_2,\ldots, x_{\ell})\bigcap I_{n+1} = (x_1,x_2,\ldots, x_{\ell})I_n$$  for all large $n.$
  Hence  $(x_1,x_2,\ldots, x_{\ell})$ is a  minimal reduction of  $\Cal F.$ The proof of (v): Let
 $\frak I$ be  a  minimal reduction  of $\Cal F.$ Now,  note that if  $x\in \frak I$ is  a weak-(FC)-element   with 
respect to
$( J, \Cal F),$  then  by (**) we have $\dim F_J(\Cal F/(x))  = \dim F_J(\Cal F) - 1.$ Consequently,  $\frak I(A/x)$ is  also  a minimal  reduction of  $\Cal F/(x).$  Hence by Proposition 2.3  and  by  induction on $\dim F_J(\Cal F)$, we easily give that     
there exists  a  maximal  weak-(FC)-sequence  $x_1,x_2,\ldots, x_{\ell}$ in $ \frak I $ with respect to $(J, \Cal F).$ By (iv), $(x_1,x_2,\ldots, x_{\ell})$ is a reduction of $\Cal F.$ Since  $(x_1,x_2,\ldots, x_{\ell}) \subset \frak I$ and 
$\frak I$ is  a  minimal reduction  of $\Cal F,$  $\frak I$ = $(x_1,x_2,\ldots, x_{\ell}).$
Proposition has been proved. 
$\square$  
\enddemo

\vskip 0.3cm 

{\bf \heading{3. The multiplicity of fiber cones of good  filtrations } 
\endheading \smallskip}
\vskip 0.2cm  
\noindent
In this section, we will   examine  the multiplicity of fiber cones of good filtrations.

Denote by $\frak M$ the maximal homogeneous ideal  of $R(\Cal F).$ Recall that  the multiplicity of the fiber  cone $F_J(\Cal F)$  is   the  Hilbert-Samuel   multiplicity  of  local ring $ F_J(\Cal F)_{\frak M}.$

We begin by the  following note.

\noindent
{\bf Note 1:}
By  Proposition 2.4,  a reduction  of $\Cal F$ is  a  minimal reduction of $\Cal F$ if  and only if  it is generated by a  weak-(FC)-sequence of the length $\ell = \ell(\Cal F).$ Let  $x_1,x_2,\ldots, x_{\ell}$ be  a   
weak-(FC)-sequence in $I_1$ 
with respect to $( J, \Cal F).$  Then   
 $\frak I =(x_1,x_2,\ldots, x_{\ell})$  is a  minimal reduction of $\Cal F.$ And  by  the proof of Proposition 2.4, $\dim F_J(\Cal F/(x_1,x_2,\ldots, x_s)) = \dim F_J(\Cal F)-s$  for  all  $s \leqslant \ell.$ This  also  means   $\frak I[A/(x_1,x_2,\ldots, x_s)]$ is a minimal reduction of  $\Cal F/(x_1,x_2,\ldots, x_s).$

The following proposition plays  an important role in the proofs of this paper.

\proclaim{\bf Proposition 3.1}  Let $J$ be an ${\frak m}$-primary ideal of  $A$.  Let  $\Cal F = \{I_n\}_{n \ge 0}$  be a good  filtration of ideals in $A$ with 
$\ell(\Cal F) = \ell > 0.$ Let $x_1,x_2,\ldots, x_{\ell}$ be a  
weak-(FC)-sequence in $I_1$ 
with respect to $( J,  \Cal F).$  Set   
$ \frak I   =   ( x_1,x_2,\ldots,x_{\ell})$ and  $I_{-1}= 0.$ 
For any $i < \ell,$ set 
$P(i) = (x_1,\ldots,x_i) : {I_1}^{\infty},\;   \Cal F(i)  = \Cal F/(x_1,\ldots,x_i), \; \Cal F'(i)  = \Cal F /P(i).$  
 Then
\noindent 
\item{ }{\rm(i)}  
 $e(F_J(\Cal F)) = e(F_J(\Cal F(i))) = e(F_J(\Cal F'(i))).$ 

\noindent
\item{ }{\rm(ii)}  
$l_A\big[{{F_J(\Cal F'(i))}\over{\frak It F_J(\Cal F'(i))}}\big] 
\leqslant l_A\big[{{F_J(\Cal F(i))}\over {\frak It F_J(\Cal F(i))}}\big]  
\leqslant l_A\big[{{F_J(\Cal F)}\over{\frak It  F_J(\Cal F)}}\big].$ 

\noindent
\item{ }{\rm(iii)}  
$l_A\big[{{F_J(\Cal F'(i))}\over{\frak It F_J(\Cal F'(i))}}\big] 
=  l_A\big[{{F_J(\Cal F)}\over{\frak It  F_J(\Cal F)}}\big]$ 
 if  and  only  if 
$I_n \bigcap P(i) \subseteq \frak I I_{n-1 } \; (\text{\rm mod }J I_n)  
\text {  for all } 0 \leqslant n \leqslant r_{\frak I}(\Cal F).$  
\endproclaim
\vskip0.2cm
\demo{Proof}  Set 
 $ A^* =A/ 0:{I_1}^{\infty}$ and  $\Cal F^* = \Cal F/0:{I_1}^{\infty} = \{I_n^* = I_nA^*\}_{n \ge 0}.$ 
 By Remark 2.1 (ii),
$(0:{I_1}^{\infty}) \bigcap I_n = 0$ for all large enough  $n.$  This  gives    
$l_A({{{I}_n^*}\over{J{I}_n^*}}) = 
l_A\big[{{I_n }\over{J I_n + (0:{I_1}^{\infty})\bigcap I_n}}\big] =
 l_A({{I_n}\over{JI_n}})$ for all large enough $n.$ Hence
$$e(F_J(\Cal F)) = e(F_J (\Cal F^*)). \eqno(2)$$
Let $x\in I_1$ be  a weak-(FC)-element  with 
respect to $( J, \Cal F).$ Then on the one  hand  by Proposition 2.4(ii),  $e(F_J(\Cal F/(x))) = e(F_J(\Cal F)).$   
Now, assume that 
the analytic spread $\ell = \ell(\Cal F) > 1$ and $x_1,x_2,\ldots,x_i$ 
$(i < \ell)$ is a weak-(FC)-sequence
in $I_1$ with respect to $( J, \Cal F).$
 We  easily  show  by induction on $i < \ell$  that
$e(F_J(\Cal F)) = e(F_J(\Cal F(i))).$
On the other hand by (2), $ e(F_J(\Cal F(i))) = e(F_J(\Cal F'(i))).$  
Hence $$e(F_J(\Cal F)) = e(F_J(\Cal F(i))) = e(F_J(\Cal F'(i))).$$
 Note that this equation is 
true too in the case of $\ell = 1$ by (2). This establishes (i). 
Set 
$r_{\frak I}(\Cal F) = r,  $ 
 $Q(i) = (x_1,\ldots,x_i), $    
 $P(i) = Q(i) : {I_1}^{\infty}.$  Since $$ F_J(\Cal F(i)) =  
\bigoplus_{n\ge 0}\big[{{I_n}\over {Q(i) \bigcap I_n + JI_n}}\big]t^n 
\; \text {and}\; F_J(\Cal F'(i)) = 
 \bigoplus_{n\ge 0}\big[{{I_n}\over{P(i) \bigcap I_n + JI_n}}]t^n, $$
$$ \aligned\frak ItF_J(\Cal F)& =  \bigoplus_{n\ge 1}\big[{{\frak II_{n-1}+JI_n}\over{ JI_n}}]t^n\\
 \frak It F_J(\Cal F(i))& = \bigoplus_{n\ge 1}\big[{{Q(i) \bigcap I_n +\frak II_{n-1}+ JI_n}\over {Q(i) \bigcap I_n + JI_n}}]t^n\\ 
\frak It F_J(\Cal F'(i))& = \bigoplus_{n\ge 1}\big[{{P(i) \bigcap I_n +\frak II_{n-1} + JI_n}\over {P(i) \bigcap I_n + JI_n}}]t^n.\endaligned$$   
 Hence 
$$\aligned l_A\big[F_J(\Cal F)/{\frak It}F_J(\Cal F)\big] &= l_A\big[A/J\big] + 
 \underset{1 \leqslant n \leqslant r}\to\sum l_A\big[{{I_n}\over{({\frak I}I_{n-1}+JI_n)}}\big]\\ 
 l_A\big[F_J(\Cal F(i))/{\frak It}F_J(\Cal F(i))\big] &= l_A\big[A/J+Q(i)\big] + 
 \underset{1 \leqslant n \leqslant r}\to\sum l_A\big[{{I_n}\over {(Q(i) \bigcap I_n +{\frak I}I_{n-1}+JI_n)}}\big]\\   
 l_A\big[F_J(\Cal F'(i))/{\frak It}F_J(\Cal F'(i))\big] &= l_A\big[A/ J+P(i)\big] +
 \underset{1 \leqslant n \leqslant r}\to\sum l_A\big[{{ I_n}\over {(P(i)\bigcap I_n + {\frak I}I_{n-1}+JI_n)}}\big]  
.\endaligned$$   
It is clear that  $$  ({\frak I}I_{n-1}+ JI_n) \subseteq 
 (Q(i)\bigcap I_n+ {\frak I}I_{n-1}+
JI_n) \subseteq   
(P(i)\bigcap I_n+{\frak I}I_{n-1}+JI_n)$$ for all $0 \leqslant n \leqslant r .$  
Hence we immediately get (ii). Moreover,    
$$ l_A\big[F_J(\Cal F)/{\frak It}F_J(\Cal F)\big] = l_A\big[F_J(\Cal F'(i))/{\frak It}F_J(\Cal F'(i))\big]$$ if and only if 
$ ({\frak I}I_{n-1}+ JI_n) = 
 (P(i)\bigcap I_n+{\frak I}I_{n-1}+JI_n) \; \text {for   all}\;  0 \leqslant n \leqslant r .$  
 This  means    
$$I_n \bigcap P(i) \subseteq \frak I I_{n-1}\; (\text{\rm mod }J I_n)  
\text {  for all } 0 \leqslant n \leqslant r.$$  
  We   have  (iii).   
$\square$ 
\enddemo

\proclaim{\bf Lemma 3.2} Let $J$ be an ${\frak m}$-primary ideal of  $A$.  Let  $\Cal F = \{I_n\}_{n \ge 0}$  be a good  filtration of ideals in $A$ with 
$\ell(\Cal F) = \ell =1$   and 
$x \in I_1$ such that $(x)$ is a reduction of $\Cal F.$  Set $r_{(x)}(\Cal F) = r.$ Then 

\noindent 
\item{ }{\rm(i)}  If  \text{\rm grade }$I_1 = 1$  then
$e(F_J(\Cal F)) = l_A(I_n/JI_n)$  for   all $n \ge r.$ 
      
\noindent
\item{ }{\rm(ii)} $e(F_J(\Cal F)) = 
l_A\big[{{I_n}\over {(0:{I_1}^{\infty})\bigcap I_n + JI_n}}\big]$ for   all $n \ge r.$  
\endproclaim

\vskip0.2cm
\demo{Proof} By Remark 2.1 (v), $x$ is  a weak-(FC)-element  in $I_1$ 
with respect to $( J,  \Cal F).$   
Since $\ell(\Cal F) = 1,$ $l_A(I_n/JI_n)$ takes a constant value 
for all large enough $n.$  This gives 
  $e(F_J(\Cal F)) = l_A(I_n/JI_n)$ for all large enough $n.$ 
Remember that $r_{(x)}(\Cal F) = r.$ 
Hence $I_n = I_rx^{n-r}$ for  all  $n > r.$ Now, if grade $I_1  > 0$ then $x^{n-r}$
is non-zero-divisor in $A.$ This  implies  the following isomorphism  of  $A$-modules $$I_r/JI_r  \simeq x^{n-r}I_r/x^{n-r}JI_r  = I_n/JI_n$$ for   all  $n \ge r.$  
We get (i).  
Set 
 $ \Cal F^* = \Cal F/0: {I_1}^{\infty} = \{I_n^* = I_n[A/0: {I_1}^{\infty}]\}_{n \ge 0}.$   Recall that by  
Proposition 3.1,  
$e(F_J(\Cal F)) = e(F_J(\Cal F^*)).$ 
Since $(0: {I_1}^{\infty}): I_1 = 0:{I_1}^{\infty},$ it follows that
grade ${I_1^*} > 0.$ 
 On the other hand  we always have  
 grade ${I_1^*} \leqslant \ell({I_1^*}) \leqslant \ell(I_1) =1.$ Consequently,     
grade ${I_1^*} = \ell({I_1^*})  =1.$
Since $(x^*) =  x[A/(0: {I_1}^{\infty})]$  is a reduction of  $\Cal F^*$ and  $r_{(x^*)}({\Cal F^*})\leqslant r_{(x)}(\Cal F)= r,$ by (i)  we  get
$$e(F_J(\Cal F^*)) = l_A(I_n^*/JI_n^*) = l_A\big[{{I_n}+(0:{I_1}^{\infty})\over {(0:{I_1}^{\infty}) + JI_n}}\big] =   
l_A\big[{{I_n}\over {(0:{I_1}^{\infty})\bigcap I_n + JI_n}}\big]$$  for   all $n \ge r.$     
Thus, $$e(F_J(\Cal F)) = l_A\big[{{I_n}\over {(0:{I_1}^{\infty})\bigcap I_n + JI_n}}\big]$$  for   all $n \ge r.$  
$\square$ 
\enddemo

 By combining Proposition 3.1 with Lemma 3.2, we obtain
the following theorem. That is the  main result of this  section.

\proclaim{\bf Theorem  3.3}  Let $J$ be an ${\frak m}$-primary ideal of  $A$.  Let  $\Cal F = \{I_n\}_{n \ge 0}$  be a good  filtration of ideals in $A$ with 
$\ell(\Cal F) = \ell > 0.$  Let $x_1,x_2,\ldots, x_{\ell}$ be a  
weak-(FC)-sequence in $I_1$ 
with respect to $( J,  \Cal F).$  Set  $$\frak I = (x_1,x_2,\ldots, x_{\ell}), \; r_{\frak I}(\Cal F) = r, \;   
  Q = (x_1,x_2,\ldots,x_{\ell-1}).$$ 
Then
 $e(F_J(\Cal F)) = 
l_A\big[{{I_n}\over {
(Q:{I_1}^{\infty})\bigcap I_n + JI_n}}\big] \text { for \; all} \; \; n \ge r.$  

\endproclaim
\vskip0.2cm

\demo{Proof} Set $\Cal F({\ell-1}) = \Cal F/Q = \{I(\ell-1)_n = I_n (A/Q) \}_{n \ge 0}.$  
     By Proposition 3.1, we get
$$e(F_J(\Cal F)) = e(F_J(\Cal F({\ell-1}))).$$  By Note 1, $\frak I $  
  is a minimal reduction of $\Cal F$ and    
$\dim F_J(\Cal F({\ell-1})) = 1,$  and  $(\bar x_{\ell}) =  x_{\ell}(A/Q) = \frak I(A/Q)$  is a minimal reduction of $\Cal F ({\ell-1}).$
Since $r_{\frak I}(\Cal F) = r,$   it follows that $r_{(\bar x_\ell)}(\Cal F(\ell-1)) \leqslant  r.$ 
 Hence by Lemma 3.2(ii),   
$$e(F_J(\Cal F(\ell-1))) = 
l_A\big[{{I(\ell-1)_n}\over {(0:{I({\ell-1})_1}^{\infty})\bigcap I(\ell-1)_n + JI(\ell-1)_n}}\big] \; \text { for   all } n \ge r.$$
Note that  we  always  have the following isomorphism  of  $A$-modules  
$${{I(\ell-1)_n}\over {(0:{I(\ell-1)_1}^{\infty})\bigcap I(\ell-1)_n + JI(\ell-1)_n}}
\simeq {{I_n + Q} \over {(Q:{I_1}^{\infty})\bigcap I_n + JI_n + Q}} \simeq {{I_n}\over {(Q:{I_1}^{\infty})\bigcap I_n + JI_n }}.$$
This gives   $$l_A\big[{{I(\ell-1)_n}\over {(0:{I(\ell-1)_1}^{\infty})\bigcap I(\ell-1)_n + JI(\ell-1)_n}}\big] =
l_A\big[{{I_n}\over {(Q:{I_1}^{\infty})\bigcap I_n + JI_n }}] .$$
Hence  $$e(F_J(\Cal F)) = 
l_A\big[{{I_n}\over {
(Q:{I_1}^{\infty})\bigcap I_n + JI_n}}\big] \text { for \; all} \;  n \ge r.$$  
Theorem 3.3 has been proved. $\square$       
\enddemo

\vskip 0.3cm 

{\bf \heading{4. The Cohen-Macaulayness of  fiber cones  of  good  filtrations}
\endheading \smallskip}

\vskip 0.2cm
\noindent
In this section,  we  answer  the question when the fiber cones of good  filtrations of ideals  in $A$  are  Cohen-Macaulay.  

Denote by $\frak M$ the maximal homogeneous ideal  of $R(\Cal F).$ Then
 $F_J(\Cal F)$ is Cohen-Macaulay if and only if $ F_J(\Cal F)_{\frak M}$
is Cohen-Macaulay [H-R].  

We begin by establishing the following lemma.

\proclaim{\bf Lemma 4.1}   Let  $\Cal F = \{I_n\}_{n \ge 0}$  be a good  filtration of ideals in $A$ with 
 with 
$\ell(\Cal F) =$ grade $I_1  = 1$ and
$x \in I_1$ such that $(x)$ is a reduction of $\Cal F.$   Set $r_{(x)}(\Cal F) = r.$ Let $J$ be an ${\frak m}$-primary ideal of  $A.$  
Then    
$F_J(\Cal F)$ is Cohen-Macaulay if and only if     
$xI_{n-1}\bigcap JI_n  = JxI_{n-1} \;
\text {for  all }\; 1 \leqslant n\leqslant r.$  
\endproclaim

\vskip0.2cm
\demo{Proof} Without loss of generality we may assume that $F_J(\Cal F) = F_J(\Cal F)_{\frak M}.$  Recall that  by Remark 2.1 (v), $x$ is  a weak-(FC)-element  in $I_1$ 
with respect to $( J,  \Cal F).$   On the one hand  by Proposition 2.4(i) and $(x)$ is a reduction of  $\Cal F,$  
$$e(F_J(\Cal F)) = e(F_J(\Cal F)^+; F_J(\Cal F)) = e(xtF_J(\Cal F); F_J(\Cal F)).$$  On the other hand  by Lemma 3.2(i),  $ e(F_J(\Cal F)) = l_A(I_r/JI_r).$ Hence $$l_A(I_r/JI_r) =  e(F_J(\Cal F)) = e(xtF_J(\Cal F); F_J(\Cal F)).$$    
Since   grade $I_1  > 0,$  $x$
is non-zero-divisor in $A.$ Hence $xI_s/xJI_s  \simeq I_s/JI_s$ for all $s \ge 0.$  Therefore,  
$$\aligned  l_A\big[F_J(\Cal F)/xtF_J(\Cal F)\big] &= 
l_A(A/J)+ \underset{1 \leqslant n \leqslant r}\to\sum l_A\big[I_n/xI_{n-1}+JI_n\big]\\ 
&= l_A(A/J)+ 
\underset{1 \leqslant n \leqslant r}\to\sum (l_A\big[I_n/JI_n\big]
-l_A\big[(xI_{n-1}+JI_n)/JI_n\big])\\ 
&= l_A(A/J)+ 
\underset{1 \leqslant n \leqslant r}\to\sum (l_A\big[I_n/JI_n\big]
-l_A\big[xI_{n-1}/xI_{n-1} \bigcap JI_n\big])\\ 
&\ge l_A(A/J)+ 
\underset{1 \leqslant n \leqslant r}\to\sum (l_A\big[I_n/JI_n\big]
-l_A\big[xI_{n-1}/xJI_{n-1}\big])\\ 
&= l_A(A/J)+ \underset{1 \leqslant n \leqslant r}\to\sum (l_A\big[I_n/JI_n\big]-l_A\big[I_{n-1}/JI_{n-1}\big]) \text{( by x is non-zero-divisor  in A)}\\ 
&= l_A\big[I_r/JI_r\big] = e(F_J(\Cal F)) = e(xtF_J(\Cal F); F_J(\Cal F)).\endaligned$$
It is clear that $xtF_J(\Cal F)$ is an  ideal  of  parameter  for  $F_J(\Cal F).$ Consequently,  $F_J(\Cal F)$  is Cohen-Macaulay 
if  and   only  if $l_A\big[F_J(\Cal F)/xtF_J(\Cal F)\big] = e(xtF_J(\Cal F); F_J(\Cal F)).$
This  is equivalent to 
$$ \underset{1 \leqslant n \leqslant r}\to\sum (l_A\big[I_n/JI_n\big]
-l_A\big[xI_{n-1}/(xI_{n-1} \bigcap JI_n)\big]) = \underset{1 \leqslant n \leqslant r}\to\sum (l_A\big[I_n/JI_n\big]
-l_A\big[xI_{n-1}/xJI_{n-1}\big]).$$ Hence      
$F_J(\Cal F)$  is Cohen-Macaulay 
if  and   only  if   
$xI_{n-1} \bigcap JI_n 
= JxI_{n-1} 
\text { for all }1 \leqslant n \leqslant r.$
$\square$ 
\enddemo

Let $J$ be an ${\frak m}$-primary ideal of  $A$ and  $\Cal F = \{I_n\}_{n \ge 0}$   a good  filtration of ideals in $A$ with 
$\ell(\Cal F) = \ell > 0.$  Remember that   if  $\ell = \ell(\Cal F)$ and  $x_1,x_2,\ldots, x_{\ell}$ is a   
weak-(FC)-sequence in $I_1$ 
with respect to $( J,  \Cal F),$ and  set $\frak I = (x_1,x_2,\ldots, x_{\ell}),$ then $\frak I$ a  minimal reduction of $\Cal F$ by Note 1. Since $\dim F_J(\Cal F) = \ell$ and $\frak ItF_J(\Cal F)$ is an $({\frak m}/J)\bigoplus F_J(\Cal F)^+$-primary ideal of  $F_J(\Cal F),$ $\frak It F_J(\Cal F)$ is  an  ideal  of  parameter  for  $F_J(\Cal F).$   

The main result of this section  is established in the following theorem. 

\proclaim{\bf Theorem  4.2} Let $J$ be an ${\frak m}$-primary ideal of  $A$.  Let  $\Cal F = \{I_n\}_{n \ge 0}$  be a good  filtration of ideals in $A$ with 
$\ell(\Cal F) = \ell > 0.$  Let $x_1,x_2,\ldots, x_{\ell}$ be a  
weak-(FC)-sequence in $I_1$ 
with respect to $( J, \Cal F).$  Set  $$\frak I = (x_1,x_2,\ldots, x_{\ell}), \; r_{\frak I}(\Cal F) = r, \;   
  Q = (x_1,x_2,\ldots,x_{\ell-1}), \; I_{-1} = 0.$$  
    Then
$F_J(\Cal F)$ is Cohen-Macaulay if and only if the following conditions are 
satisfied: 

\noindent
\item{ }{\rm(i)}     
$(Q:{I_1}^{\infty})\bigcap I_n \subseteq \frak II_{n-1}\; (\text{\rm mod }{J}I_n)$  
for all \; $0 \leqslant n\leqslant r.$

\noindent
\item{ }{\rm(ii)}     
$\big[\frak II_{n-1}+ (Q:{I_1}^{\infty})\big]\bigcap {J}I_n = J{\frak I}I_{n-1 }\;\; (\text {\rm mod }Q:{I_1}^{\infty})$ 
  for  all \; $1 \leqslant n\leqslant r.$
\endproclaim
\vskip 0.2cm
\demo{Proof}Denote by $\frak M$ the maximal homogeneous ideal  of $R(\Cal F).$ Without loss of generality we may assume that $F_J(\Cal F) = F_J(\Cal F)_{\frak M}.$ Set $\Cal F({\ell-1}) = \Cal F/Q = \{I(\ell-1)_n = I_n (A/Q) \}_{n \ge 0}$ and  
$$\Cal F'({\ell-1}) = \Cal F/Q:{I_1}^{\infty} = \{I'(\ell-1)_n = I_n (A/Q:{I_1}^{\infty}) \}_{n \ge 0}.$$  Since $\frak ItF_J(\Cal F)$ is  a  reduction of  $F_J(\Cal F)^+,$ $\frak It F_J(\Cal F'({\ell-1}))$ is a  reduction of   $F_J(\Cal F'({\ell-1}))^+. $ Hence 
by Proposition 2.4(i) and  Proposition 3.1,  we have
$$\aligned  e(\frak It F_J(\Cal F); F_J(\Cal F)) =   e(F_J(\Cal F))  = e(F_J(\Cal F'({\ell-1})))& = e(\frak ItF_J(\Cal F'({\ell-1})); F_J(\Cal F'({\ell-1})))\;\; \text {and }\\ 
 l_A\big[F_J(\Cal F'({\ell-1}))/\frak ItF_J(\Cal F'({\ell-1}))\big] &\leqslant
 l_A\big[F_J(\Cal F)/\frak ItF_J(\Cal F)\big].\endaligned \eqno (3)$$ 
Since $(Q: {I_1}^{\infty}): I_1 = Q:{I_1}^{\infty},$ it follows that
 grade $I'({\ell-1})_1  > 0.$ By Note 1, $\ell(I({\ell-1})_1) = 1.$  Hence $\ell(I'({\ell-1})_1)\leqslant \ell(I({\ell-1})_1) = 1.$ On the other hand,    
 grade $I'({\ell-1})_1 \leqslant \ell(I'({\ell-1})_1)$  is  always  true.   Consequently,   
grade $I'({\ell-1})_1 = \ell(I'({\ell-1})_1)  = 1.$  
 Since   $x_{\ell}tF_J(\Cal F'({\ell-1})) =  \frak It F_J(\Cal F'({\ell-1}))$ is  a  reduction of  $F_J(\Cal F'({\ell-1}))^+$ and 
$\ell(I'({\ell-1})_1)  = 1,$
it follows that $\frak It F_J(\Cal F'({\ell-1}))$ is  an  ideal  of  parameter  for  $F_J(\Cal F'({\ell-1})).$ This  gives   
$$ e(\frak ItF_J(\Cal F'({\ell-1})); F_J(\Cal F'({\ell-1}))) \leqslant  l_A\big[F_J(\Cal F'({\ell-1}))/\frak ItF_J(\Cal F'({\ell-1}))  \big].\eqno (4)$$
   Since  $\frak ItF_J(\Cal F)$ is  an  ideal  of  parameter  for  $F_J(\Cal F),$   $F_J(\Cal F)$  
is Cohen-Macaulay if and only if    
$$ e(\frak ItF_J(\Cal F); F_J(\Cal F)) = l_A\big[F_J(\Cal F)/\frak ItF_J(\Cal F)\big].\eqno (5)$$ 
By (3) and (4),  (5)  is equivalent to  the  following two equations:
 $$ e(\frak ItF_J(\Cal F'({\ell-1})); F_J(\Cal F'({\ell-1}))) = l_A\big[F_J(\Cal F'({\ell-1}))/\frak ItF_J(\Cal F'({\ell-1}))\big]  \eqno (6) $$  and  
 $$l_A\big[F_J(\Cal F'({\ell-1}))/\frak ItF_J(\Cal F'({\ell-1}))\big]  = l_A\big[F_J(\Cal F)/\frak ItF_J(\Cal F)\big].\eqno (7)$$
Recall  that  $\frak ItF_J(\Cal F'({\ell-1}))$ is  an  ideal  of  parameter  for  $F_J(\Cal F'({\ell-1})).$  Hence (6) is  equivalent to
$F_J(\Cal F'({\ell-1}))$
  is Cohen-Macaulay. 
On the one hand,     
$$ \text {grade} I'({\ell-1})_1 = \ell(I'({\ell-1})_1)  = 1 \;
  \text {and }\;     
x_{\ell}tF_J(\Cal F'({\ell-1})) = \frak It F_J(\Cal F'({\ell-1}))$$ is a reduction of  $F_J(\Cal F'({\ell-1}))^+.$ On the other  hand, $r_{x_{\ell}(A/Q:{I_1}^{\infty})}(\Cal F'({\ell-1}))  \leqslant r_{\frak I}(\Cal F) =r.$  Hence by Lemma 4.1, 
$F_J(\Cal F'({\ell-1}))$  is Cohen-Macaulay if and only if  $$x_{\ell}I'({\ell-1})_{n-1}\bigcap JI'({\ell-1})_{n} = x_{\ell}JI'({\ell-1})_{n-1} \;
\text {for  all }\; 1 \leqslant n\leqslant r.\eqno (8)$$  Since $x_{\ell}tF_J(\Cal F'({\ell-1})) = \frak ItF_J(\Cal F'({\ell-1}))$, (8)  means    
$$\big[\frak II_{n-1}+ (Q:{I_1}^{\infty})]\bigcap \big[(Q:{I_1}^{\infty})+{J}I_n\big] = J{\frak I}I_{n-1 }+ (Q:{I_1}^{\infty})\; 
\text {for  all}\; 1 \leqslant n\leqslant r.$$
It can be verified that 
this  condition  also  means   
$$\big[\frak II_{n-1}+ (Q:{I_1}^{\infty})]\bigcap {J}I_n = J{\frak I}I_{n-1 } (\text {\rm mod }Q:{I_1}^{\infty})\; \text{for  all}\; 1 \leqslant n\leqslant r. $$ 
By Proposition  3.1(iii), (7) is equivalent to   
$$I_n \bigcap (Q : {I_1}^{\infty}) \subseteq \frak I I_{n-1} \;  (\text{\rm mod }J I_n)  
\text {  for all } 0 \leqslant n \leqslant r.$$
Hence
$F_J(\Cal F)$ is Cohen-Macaulay if and only if the following conditions 
are satisfied
$$\aligned    
(Q:{I_1}^{\infty})\bigcap I_n &\subseteq \frak II_{n-1}\; (\text{\rm mod }{J}I_n)  
\; \text {for  all } 0 \leqslant n\leqslant r .\\
\big[\frak II_{n-1}+ (Q:{I_1}^{\infty})\big]\bigcap {J}I_n &= J{\frak I}I_{n-1 } (\text {\rm mod }Q:{I_1}^{\infty}) \;
\text {for   all }\; 1 \leqslant n\leqslant r. \endaligned$$ 
This completes the proof of Theorem  4.2. $\square$ 

     From  Theorem  4.2, we give the following  interesting 
consequence.  

\proclaim{\bf Corollary 4.3} Let $(A, \frak m)$ be a Noetherian local ring 
with  maximal ideal $\frak m,$  $\Cal  F = \{I_n\}_{n \ge 0}$  a good   filtration of ideals in $A$  with 
 $\ell(\Cal F) =  1.$  Set $\Cal  F^{(T)} = \{I_{Tn} \}_{n \ge 0}.$
 Then    
$F_{\frak m}(\Cal  F^{(T)})$ is  Cohen-Macaulay for  all   large \;$T.$     
\endproclaim

\vskip0.2cm
\demo{Proof} Set $A^* = A/ 0:{I_1}^{\infty}, {\frak m^*}={\frak m}A^*, \Cal  F^* = \Cal F/0:{I_1}^{\infty} = 
\{I_n^* = I_nA^*\}_{n \ge 0}.$ Let $x$ be an   element in  $I_1$  such that $(x)$ is  a  minimal  reduction  of $\Cal F.$ Then  $(x^*)$ ($x^*$  the image of $x$ in $A^*$) is  also  a  minimal  reduction  of $\Cal F^*.$ Hence by Remark 2.1(v),  $x$ and  $x^*$  are  
weak-(FC)-elements  with respect
to $({\frak m}, \Cal F)$ and  $({\frak m^*}, \Cal F^*),$ respectively.   Since $x^*$ satisfies the condition (FC$_1$),  there exists  a positive  integer  $u$  such that  
$(x^*)\bigcap {\frak m^*}^mI_n^*  = {\frak m^*}^mx^*I_{n-1}^*$ for  all $n \ge u$      
 and  all non-negative integers $m.$ By Remark 2.1 (ii),
there exists  a positive  integer  $v$ such that 
$(0: {I_1}^{\infty}) \bigcap I_n  = 0$ for  all $n \ge v.$
Set $N = \max\{u, v \}.$  
 For  the  proof of this  corollary  we need  the following.  

\noindent
{\bf Note 2:}
  Let  $C$ and $D$ be  subsets of  $A^*.$ 
   If  $y \in {x^*} C \bigcap x^*D$  then $ y = x^*a = x^*b$  with $a \in C ; b \in D.$  Hence $x^*(a-b) = 0.$   This means 
 $(a-b) \in (0: x^*).$ Since $x$ satisfies the condition (FC$_2$),  $x^*$ is  a non-zero-divisor in $A^*.$   
Consequently,
 $(a-b) \in (0: x^*)  = 0.$   
Hence $a =b \in C\bigcap D.$ This  implies  that $y \in  x^*(C \bigcap D).$ Thus, $x^* C \bigcap x^* D = x^*(C \bigcap D).$

  Now we  choose $T \ge N.$  Using induction on $h \ge 1 ,$ we will prove that   
$${x^*}^hI_{T(n-1)}^*\bigcap {\frak m^*}^mI_{T(n -1)+h}^* = {\frak m^*}^m{x^*}^hI_{T(n-1)}^*\eqno (9)$$  for  all $n \ge 2$ and     
 for all non-negative integers $m.$ Since $n \ge 2$,  $T(n-1) \ge T \ge N$. Therefore, 
$$x^*I_{T(n-1)}^*\bigcap {\frak m^*}^mI_{T(n -1)+1}^* \subseteq (x^*)\bigcap {\frak m^*}^mI_{T(n -1)+1}^* =  
{\frak m^*}^mx^*I_{T(n-1)}.$$ Hence $$x^*I_{T(n-1)}^*\bigcap {\frak m^*}^mI_{T(n -1)+1}^* =  {\frak m^*}^mx^*I_{T(n-1)}^*.\eqno (10)$$ 
Thus, (9)  is true for $h=1.$ Suppose that (9)  was  true for  $h-1 \ge 1.$ Since  obvious  facts 
$${x^*}^hI_{T(n-1)}^* \subset x^*I_{T(n-1)+h-1}^* \subset I_{T(n-1)+h}^*\; \text { and }\; T(n-1)+h - 1 > N,$$  
$$ \aligned {x^*}^hI_{T(n-1)}^* \bigcap {\frak m^*}^mI_{T(n-1)+h}^*  &= 
{x^*}^hI_{T(n-1)}^*\bigcap \big[x^*I_{T(n-1)+h-1}^* \bigcap {\frak m^*}^mI_{T(n-1)+h}^*\big]\\
& = {x^*}^hI_{T(n-1)}^*\bigcap  {\frak m^*}^mx^*I_{T(n-1)+h-1}^* \text { (by}\;\; (10))\\
 & = x^*\big[{x^*}^{h-1}I_{T(n-1)}^*\bigcap  {\frak m^*}^mI_{T(n-1)+h-1}^*\big]\text {(by Note 2)}. \endaligned$$
By (9) is true for $h-1,$  $${x^*}^{h-1}I_{T(n-1)}^*\bigcap  {\frak m^*}^mI_{T(n-1)+h-1}^* =  {\frak m^*}^m{x^*}^{h-1}I_{T(n-1)}^*.$$ Hence
${x^*}^hI_{T(n-1)}^*\bigcap {\frak m^*}^mI_{T(n-1)+h}^* =  {\frak m^*}^m{x^*}^hI_{T(n-1)}^*.$  Consequently, the induction is complete.  So $${x^*}^hI_{T(n-1)}^*\bigcap {\frak m^*}^mI_{T(n -1)+h}^* =  {\frak m^*}^m{x^*}^hI_{T(n-1)}^*$$  for  all $n \ge 2$ and   for all non-negative integers $m.$ This  gives that for  any $T \ge N,$    
$${x^*}^TI_{T(n-1)}^*\bigcap {\frak m^*}^mI_{Tn}^* =  {\frak m^*}^m{x^*}^TI_{T(n-1)}^*$$  for  all $n \ge 2$  and   for all non-negative integers $m.$ Hence                                                                  
${x^*}^TI_{T(n-1)}^*\bigcap {\frak m^*}I_{Tn}^* =  {\frak m^*}{x^*}^TI_{T(n-1)}^*$(***)  for  all $n \ge 2.$ Note that $(x^*)$ is a minimal reduction of $\Cal F^*,$ 
 $({x^*}^T)$ is  a  minimal  reduction of $\{I_{Tn}^* \}_{n \ge 0}.$ Now, we  need  to  show that     
$({x^*}^T)\bigcap {\frak m^*}I_T^* = {\frak m^*}{x^*}^T$  by 
 using the following note. 

\noindent
{\bf Note 3:} 
 Let $ (x_1,x_2,\ldots, x_{\ell})$ be a  minimal reduction of a good   filtration $\Cal  F = \{I_n\}_{n \ge 0}$.  
 Set $ \frak I   =   (x_1,x_2,\ldots,x_{\ell}).$ 
Since $\frak I$  is a minimal reduction of $\Cal F,$ there exist elements 
 $y_1,\ldots,y_s \in I_1$  such that $x_1, x_2, \ldots,x_{\ell}, y_1,\ldots,y_s \in I_1$ 
is a minimal base of  $I_1.$ 
Now assume that $x = a_1x_1+\cdots+a_{\ell}x_{\ell} \in
\frak I \bigcap {\frak m}I_1,$ $$x = b_1x_1+\cdots+b_{\ell}x_{\ell}+c_1y_1+\cdots+c_sy_s$$  
where   $b_i \in {\frak m}$ for $1 \leqslant  i \leqslant \ell$ and $c_j \in {\frak m}$
for $1 \leqslant j \leqslant s.$ Then we have
$$(a_1- b_1)x_1+\cdots+(a_{\ell}-b_{\ell})x_{\ell}+c_1y_1+\cdots+c_sy_s = 0.$$  
Since $x_1,x_2, \ldots,x_{\ell}, y_1,\ldots,y_s$ is a minimal base of $I_1,$ 
 $a_i-b_i \in {\frak m}$ for $1 \leqslant i \leqslant \ell.$ Hence
 $a_i \in {\frak m}$ for $1 \leqslant i \leqslant \ell.$ Thus,  $x \in {\frak m}\frak I.$   This  follows  that
$\frak I\bigcap {\frak m}I_1 = {\frak m}\frak I.$ 
      
On the one hand  by  Note 3, $({x^*}^T)\bigcap {\frak m^*}I_T^* = {\frak m^*}{x^*}^T.$  Hence ${x^*}^TI_{T(n-1)}^*\bigcap {\frak m^*}I_{Tn}^* =
  {\frak m^*}{x^*}^TI_{T(n-1)}^*$ is true for $n=1$.  On  the other hand  by (***), 
${x^*}^TI_{T(n-1)}^*\bigcap {\frak m^*}I_{Tn}^* =  {\frak m^*}{x^*}^TI_{T(n-1)}^* \; \text { for  all} \; n \ge 2.$
Consequently,  $${x^*}^TI_{T(n-1)}^*\bigcap {\frak m^*}I_{Tn}^* =  {\frak m^*}{x^*}^TI_{T(n-1)}^*\; \text { for  all} \; n \ge 1. \eqno (11)$$ 
 Since $x^T$ is  a minimal  reduction  of $\Cal F^{(T)}$, $x^T$  is  a weak-(FC)-element with  respect to  $(\frak m, \Cal F^{(T)}$) by Remark 2.1 (v).  Hence  by  Theorem 4.2,    
$F_{\frak m}(\Cal F^{(T)})$ is  Cohen-Macaulay if and only if  $\Cal F^{(T)}$  satisfies the following conditions 

\noindent
\item{ }{\rm(a)}     
$(0:{{I_{T}}}^{\infty})\bigcap I_{Tn} \subseteq  x^TI_{T(n-1)}\; (\text{\rm mod }{\frak m}I_{Tn})$  
for all \; $0 \leqslant n\leqslant r,$ where $r$ is  the reduction  number  of  $\Cal F^{(T)}.$

{\noindent
\item{ }{\rm(b)}     
$\big[x^TI_{T(n-1)}+ (0:{I_T}^{\infty}) \big]\bigcap {\frak m}I_{Tn}  = {\frak m}x^TI_{T(n-1)} (\text {\rm mod }0:{{I_T}}^{\infty})$ 
  for  all \; $1 \leqslant n\leqslant r.$}  
Since $Tn \ge N$ for  all $n \ge 1$ \;  and  $0: {I_T}^{\infty} = 0: {I_1}^{\infty}$,  $(0: {I_T}^{\infty}) \bigcap I_{Tn} = 0$ for  all $n \ge 1.$   Remember that in the statement of Theorem 4.2,  we  assigned  ${I_T}_{-1} = 0.$ This means $I_{T(-1)} = 0.$  Hence conditions (a) and (b) are   
equivalent  to 
$$\text {(a'):}\;     
(0:{{I_{T}}}^{\infty})\bigcap A =  0 \; (\text{mod }{\frak m})  
\; \text {and }\; {\text (b'):}\;     
{x^*}^TI_{T(n-1)}^*\bigcap {\frak m^*}I_{Tn}^*  = {\frak m^*}{x^*}^TI_{T(n-1)}^* 
  \; \text {for   all} \; 1 \leqslant n\leqslant r,$$ respectively.
But (a') is obvious and (b') is satisfied  by (11).
Thus, 
$F_{\frak m}(\Cal F^{(T)})$ is  Cohen-Macaulay for  all   large \;$T,$  as  required. $\square$
\enddemo               

\vskip 0.2 cm
\noindent {\bf Remark  4.4.} Note that in  the  case  where $\Cal F = \{I^n\}_{n \ge 0}$  is  an $I$-adic filtration.  Set $A^*=\dfrac{A}{0:I^\infty}$, $\frak m^* = \frak m  A^*$,  ${I}^* = IA^*,$ and $x^*$ the image of $x \in I $ in $A^*.$  Recall that  the  condition (FC$_1$) of  [Definition  in  Sect. 2, Vi2]  is $(x^*)\cap {\frak m^*}^{m}{I^*}^{n} =  x^* {\frak m^*}^{m}{I^*}^{n-1}$ for  all   large $n$ and  for  all  non-negative integers $m.$ So  the condition (FC$_1$) in this  paper  is a weaker condition than the condition
(FC$_1$) in [Vi2].  Hence  as   immediate consequences  of Theorem 3.3 and Theorem 4.2, we obtained   more  favorite
results  than the results in  [Vi2].

\vskip 0.4cm
\Refs
\vskip 0.2cm
\item{[C]}  T.  Cortadellas, {\it Fiber cones with almost maximal depth} , Comm. Algebra
33 (2005), no. 3, 953–963.

\item{ [C-G-P-U]}  A. Corso, L. Ghezzi, C. Polini, B. Ulrich, 
{\it  Cohen-Macaulayness of special fiber rings},
 Comm. Algebra 31(2003), 
3713-3734. 
\vskip 0.2cm
\item{ [C-P-V]}  A. Corso, C. Polini, W. V. Vasconcelos,
{\it Multiplicity of the special fiber ring of blowups},
 Math. Proc. Camb. Phil. Soc. 140(2006), 207-219. 

\vskip 0.2cm
\item{ [C-Z]}  T. Cortadellas, S. Zarzuela,  
{\it  On the depth of the fiber cone of filtrations},  
J. Algebra 198(1997), 428-445. 
\vskip 0.2cm
\item{ [Co-Z]}  T. Cortadellas, S. Zarzuela,  
{\it  On the structure of   the  fiber  cone  of  ideals  with  analytic  spread  one},  
arXiv: math.AC/0603042 v1  2 Mar  2006.  
\vskip 0.2cm
\item{ [D-G-H]} M. D'Anna, A. Guerrieri, W. Heinzer, 
{\it Ideals  having one-dimensional fiber cone},
Lecture Notes in Pure Appl. Math. 220(2001), 155-170.  
\vskip 0.2cm
\item{ [D-R-V]} C. D'Cruz, K.N. Raghavan, J.K. Verma, 
{\it Cohen-Macaulay fiber cones}, Commutative algebra, algebraic geometry and
computational methods (Hanoi, 1996), 233-246, Springer, Singapore,1999.
\vskip 0.2cm
\item{ [G-S]}  S. Goto, Y. Shimoda, 
{\it  On the Rees algebras of Cohen-Macaulay local rings}, Lect. 
Notes in Pure and Appl. Math. 68(1979), Dekker, 201-231.
\vskip 0.2cm
\item{ [H-H]} R. H\"ubl, C. Huneke, 
{\it Fiber cones and the integral closure of ideals},
Collect. Math. 52(2001),85-100. 
\vskip 0.2cm
\item{[H-R]} M. Hochter, M. Ratliff, 
{\it Five theorems on Cohen-Macaulay rings }, Pac. J. Math. 44(1973), 147-173.
\vskip 0.2cm
\item{[H-Sa]}  C. Huneke, J. D. Sally, {\it  Birational extensions 
in dimension two and integrally closed ideals},  J. Algebra 115(1988), 481-500.
\vskip 0.2cm
\item{[H-S]}  C.  Huneke,  I.  Swanson, {\it  Integral Closure of Ideals, Rings, and Modules},
London Mathematical Lecture Note Series 336, Cambridge University Press
(2006)
\vskip 0.2cm
\item{[K]}  D.  Kirby, {\it  A note on  superficial   elements of an ideal in a local ring,}  Quart.
J. Math. Oxford (2), 14 (1963), 21-28.

\vskip 0.2cm
\item{[J-V-S]} A.V. Jayanthan, J.K. Verma, B. Singh, 
{\it Hilbert  coefficients  and  depth  of  form  rings}, 
Comm. Algebra 32(2004),1445-1452. 
\vskip 0.2cm
\item{[J-V1]} A.V. Jayanthan, J.K. Verma, 
{\it Hilbert coefficients and depth of fiber cones}, 
J. Pure Appl. Algebra 201(2005),97-115. 
\vskip 0.2cm
\item{[J-V2]} A.V. Jayanthan, J.K. Verma, 
{\it Fiber cones of ideals with almost minimal multiplicity}, 
Nagoya Math. J. 177(2005),155-179. 
\vskip 0.2cm
\item{[J-P-V]} A.V. Jayanthan, T.J. Puthenpurakal, J.K. Verma, 
 {\it On fiber cones of {\rm  m}-primary ideals}, 
Canad. J. Math. 59(2007),109-120. 
\vskip 0.2cm
\item{[N-R]}  D. G. Northcott and D. Rees, {\it  Reduction of ideals in local rings},   Proc. Cambridge Phil. Soc. 50 (1954), 145-158.
\vskip 0.2cm
\item{[Re]}  D. Rees, {\it  Generalizations of reductions and mixed 
multiplicities},  J. London Math. Soc. 29(1984), 397-414. 

\vskip 0.2cm
\item{[R-V]}  M. E. Rossi, G. Valla, {\it  Hilbert Function of Filtered Modules, } (2008)
arXiv:0710.2346

\vskip 0.2cm
\item{[Sa]} J. D. Sally, {\it  On the associated graded ring of a local Cohen-Macaulay ring},  J. Math. Kyoto  17(1977), 19-21.
\vskip 0.2cm
\item{[Sh1]}  K. Shah, {\it On the Cohen-Macaulayness of the fiber cone of an ideal},  J. Algebra 143(1991), 156-172.
\vskip 0.2cm

\item{[Sh2]}  K. Shah, {\it On equimultiple ideals},  Math. Z. 215(1994), 13-24.
\vskip 0.2cm
\item{[Vi1]}  D. Q. Viet, Mixed multiplicities of arbitrary ideals in local rings, Comm. Algebra 28(8)(2000), 3803-3821.

\vskip 0.2cm
\item{ [Vi2]}  D. Q. Viet, On the Cohen-Macaulayness of   fiber cones, Proc. Amer. Math. Soc. 136 (2008), 4185-4195.
\vskip 0.2cm

\item{[Z-S]}  O. Zariski, P.  Samuel, {\it  Commutative Algebra, Vol II, Van Nostrand,}  New York,
1960.

\end{document}